\documentclass[12pt]{article}
\usepackage[english]{babel}
\usepackage[centertags]{amsmath}
\usepackage{amsfonts}
\usepackage{amssymb}
\usepackage{amsthm}
\usepackage{newlfont}
\usepackage[bookmarksnumbered, colorlinks, plainpages]{hyperref}

\usepackage{natbib}

\begin{document}
\bibliographystyle{plainnat}

\title{On  the numerical radii of $2\times2$  complex matrices}
\author{L. Z. Gevorgyan\\ State Engineering University of
Armenia,\\ Department of Mathematics, E-mail: levgev@hotmail.com}
\date{}
\maketitle

AMS Mathematics Subject Classification {47A12,\{Primary\}, 51N15, \{Secondary\}}

\begin{abstract}
The numerical radius of the general $2\times2$ complex matrix is calculated.  \end{abstract}

\section{Introduction and preliminaries}

Let  $A$ be a linear bounded operator, acting in a Hilbert
space $\left ( {\mathcal{H},\left\langle { \bullet , \bullet }
\right\rangle }\right).$ Denote by $ W\left( A \right) = \left\{ {\left\langle {Ax,x}
\right\rangle :\left\| x \right\| = 1} \right\}$ the numerical range and  by
$w\left( A \right) = \mathop {\sup }\limits_{\lambda  \in W\left( A \right)} \left| \lambda  \right|$ the numerical radius of the operator $A.$ As it is well known, the numerical radius defines a norm on the space of all bounded operators, equivalent to the usual one, i.e.\[
w\left( A \right) \leqslant \left\| A \right\| \leqslant 2w\left( A \right).\]
The importance of the numerical radius is partially motivated by its role in the
investigation of many iterative schemes of solution of operator equations \citet{MR1206415}. Despite the great abundance of research papers on the numerical range and the numerical radius, the exact calculations of the numerical radius of general non-normal operator are very scarce. A few known results concern the $n-$ dimensional Jordan block
\[
J_n  = \left( {\begin{array}{*{20}c}
   0 & 1 & 0 &  \cdots  & 0  \\
   0 & 0 & 1 &  \cdots  & 0  \\
   \cdot & \cdot & \cdot & \cdot & \cdot  \\
   0 & 0 & 0 &  \cdots  & 1  \\
   0 & 0 & 0 &  \cdots  & 0  \\

 \end{array} } \right)
\]
($w\left( {J_n } \right) = \cos \frac{\pi }
{{n + 1}}$)
and the Volterra integration operator \[
\left( {Vf} \right)\left( x \right) = \int\limits_0^x {f\left( t \right)} dt,f \in L^2 \left( {0;\,1} \right)
\]
($w\left(V\right)=1/2$).

In the sequel we consider the simplest case of two-dimensional  space.
\section{Main result}

As any operator in a finite-dimensional space may be reduced to the Schur's upper triangular form, we get  the matrix \[
A = \left( {\begin{array}{*{20}c}
   {\lambda _1 } & {\lambda _3 }  \\
   0 & {\lambda _2 }  \\

 \end{array} } \right).\]
 Recall that the numerical range of an operator, acting in a two-dimensional space is an elliptical disk (\citet{MR675952}, ch. 22) having as foci two eigenvalues $\left\{ {\lambda _1 ,\lambda _2 } \right\}$ and \[
\left\{ {tr\left( {A^* A} \right) - \left| {\lambda _1 } \right|^2  - \left| {\lambda _2 } \right|^2 } \right\}^{1/2}\] as minor axis $2b.$

We have  $b = \left| {\lambda _3 } \right|/2,$ the distance between foci is
$2c = \left| {\lambda _2  - \lambda _1 } \right|,$ hence major axis is
$2a = \sqrt {\left| {\lambda _2  - \lambda _1 } \right|^2  + \left| {\lambda _3 } \right|^2 }.$ The center of symmetry of the ellipse has the affix
$m=\frac{{\lambda _2  + \lambda _1 }}{2}$ and whole ellipse is rotated by the angle $\varphi  = \arg \left( {\lambda _2  - \lambda _1 } \right).$
The equation of the ellipse is \[
\begin{gathered}
  \frac{{\left( {\left( {x - \operatorname{Re} m} \right)\cos \phi  + \left( {y - \operatorname{Im} m} \right)\sin \phi } \right)^2 }}
{{a^2 }} +  \hfill \\
  \frac{{\left( {\left( {x - \operatorname{Re} m} \right)\sin \phi  - \left( {y - \operatorname{Im} m} \right)\cos \phi } \right)^2 }}
{{b^2 }} - 1 = 0. \hfill \\
\end{gathered}
\]

After elementary transformations it may be reduced to the form
\[G\left(x,y \right)=Ax^2  + 2Bxy + Cy^2  + 2Dx + 2Ey + F = 0,\] where
\[
\begin{gathered}
  A = \frac{{\cos ^2 \phi }}
{{a^2 }} + \frac{{\sin ^2 \phi }}
{{b^2 }},B = \sin \phi \cos \phi \left( {\frac{1}
{{a^2 }} - \frac{1}
{{b^2 }}} \right),C = \frac{{\sin ^2 \phi }}
{{a^2 }} + \frac{{\cos ^2 \phi }}
{{b^2 }}, \hfill \\
  D =  - A\operatorname{Re} m - B\operatorname{Im} m,E =  - B\operatorname{Re} m - C\operatorname{Im} m, \hfill \\
  F = \left( {\frac{{\cos \phi \operatorname{Re} m + \sin \phi \operatorname{Im} m}}
{a}} \right)^2  + \left( {\frac{{\cos \phi \operatorname{Im} m - \sin \phi \operatorname{Re} m}}
{b}} \right)^2  - 1. \hfill \\
\end{gathered}
\]

The square of numerical radius $w\left(A\right)$ may be found as the greatest value of the sum $x^2  + y^2,$ where $x$ and $y$ are constrained by the equation $G\left(x,y\right)=0.$ According to the Lagrange multiplier method, the auxiliary function \[
L\left( {x,y,\lambda} \right) = x^2  + y^2  + \lambda G\left( {x,y} \right)
\] should be introduced and  the system of equations
\[
\begin{cases}
  \frac{{\partial L}}{{\partial x}} = 0, \\
  \frac{{\partial L}}{{\partial y}} = 0, \\
  \frac{{\partial L}}{{\partial \lambda }} = 0
\end{cases}
\]
will be solved.

We prefer another approach. The equation of the circumference is \[\Gamma \left(x,y \right) =x^2  + y^2  - R^2  = 0.\]
If the ellipse lies inside the circle and touches the circumference, then they have a common tangent, therefore  \[
\frac{{Ax + By + D}}
{{Bx + Cy + E}} = \frac{x}
{y}.\]

Finally, we get the set of simultaneous equations \begin{equation}
\left\{ \begin{gathered}
  Ax^2  + 2Bxy + Cy^2  + 2Dx + 2Ey + F = 0, \hfill \\
  Bx^2  + \left( {C - A} \right)xy - By^2  + Ex - Dy = 0. \hfill \\
\end{gathered}  \right.
\label{tan}\end{equation}
Note that the second curve is a hyperbola. The solution of this set may be reduced to a quartic polynomial equation. Using some notions of projective geometry (see \citet{MR2791970}, ch. 11 ), it is possible to avoid cumbersome calculations. Introducing the homogeneous coordinates $\left( {\begin{array}{*{20}c} x & y & 1  \\ \end{array} } \right)^t$
 of  the point $\left( {x;y} \right), $ two matrices \[
G = \left( {\begin{array}{*{20}c}
   A & B & D  \\
   B & C & E  \\
   D & E & F  \\

 \end{array} } \right)\,\,{\text{and  }}
H = \frac{1}
{2}\left( {\begin{array}{*{20}c}
   {2B} & {C - A} & E  \\
   {C - A} & { - 2B} & { - D}  \\
   E & { - D} & 0  \\

 \end{array} } \right)\]
may be associated with \eqref{tan}. Two conics define a pencil of conics $\alpha H + \beta G;\alpha ,\beta  \in \mathbb{R},$  passing through the intersection points of $G$ and $H.$ Solving the generalized eigenvalue problem $\det \left( {H - \lambda G} \right) = 0$ we get the degenerated conics $C = H - \lambda _0 G,\det \left( C \right) = 0.$ Now we calculate the adjugate matrix $J=Ad\left( C \right).$ Maybe the simplest way is the use of the formula $
J =  - \left( {a_1 I + a_2 C + C^2 } \right),$ where $a_1  = \frac{1}
{2}\left( {tr^2 \left( C \right) - tr\left( {C^2 } \right)} \right),a_2  =  - tr\left( C \right).$ Denote $p_1  = \sqrt {J_{11} } ,p_2  = \sqrt {J_{22} } ,p_3  = \sqrt {J_{33} },$ where the signs before the radicals are chosen in a such way that
$\left\{ {p_1 ,p_2 ,p_3 } \right\}^T  \cdot \left\{ {p_1 ,p_2 ,p_3 } \right\} = J$
and construct the antisymmetric matrix \[
P = \left( {\begin{array}{*{20}c}
   0 & {p_3 } & { - p_2 }  \\
   { - p_3 } & 0 & {p_1 }  \\
   {p_2 } & { - p_1 } & 0  \\
 \end{array} } \right).\]

Now we choose the parameter $\alpha$ to satisfy the condition $K = C - \alpha P,rank\left( K \right) = 1.$ Any row (or column) of $K$ defines a degenerate conics (a straight line) $l.$ Replacing the first equation in \eqref{tan} by $l=0$ we may found easily (solving a quadratic equation) the coordinates of the intersection points. One of them will supply the minimum of $\sqrt {x^2  + y^2 }$ and the second point is at the distance $w\left( A \right)$ from the origin of the coordinate system. 

We append below a MatLab program, which admits as input a $2\times2$ upper triangular matrix and calculates its numerical radius.

\section{Appendix}

 \quad function numradius(A)

nor=norm(A);

b=.5*abs(A(1,2));

c=.5*abs(A(1,1)-A(2,2));

a=sqrt(b\^\,2+c\^\,2);

m=.5*(A(1,1)+A(2,2));

k=real(m);l=imag(m);

phi=angle(A(2,2)-A(1,1));

s=linspace(0,2*pi,2000);

u=a*cos(s);v=b*sin(s);

R=[cos(phi) -sin(phi);sin(phi) cos(phi)];

K=R*[u;v];

z=K(1,:)+k+i*(K(2,:)+l);

plot(real(z),imag(z),'r')

axis 'equal'

A=cos(phi)\^\,2/a\^\,2+sin(phi)\^\,2/b\^\,2;

B=sin(phi)* cos(phi)*(1/a\^\,2-1/b\^\,2);

C=cos(phi)\^\,2/b\^\,2+sin(phi)\^\,2/a\^\,2;

D=-k*A-l*B;

E=-k*B-l*C;

F=(k*cos(phi)+l*sin(phi))\^\,2/a\^\,2+(l*cos(phi)-k*sin(phi))\^\,2/b\^\,2-1;

G=[A B D;B C E;D E F];

H=.5*[2*B C-A E;C-A -2*B -D;E -D 0];

e=eig(H,G);

for j=1:3

    if isreal(e(j))

        c=e(1);

    end

end

C=H-c*G;

s1=trace(C);s2=trace(C\^\,2);a2=-s1;a1=-.5*(s1*a2+s2);

Ad=-(a1*eye(3)+a2*C+C\^\,2);

p(1)=sqrt(Ad(1,1));

p(2)=sign(Ad(1,2))*sqrt(Ad(2,2));

p(3)=sign(Ad(1,3))*sqrt(Ad(3,3));

P=[0 p(3) -p(2);-p(3) 0 p(1);p(2) -p(1) 0];

e=eig(C,P);

e(4)=sqrt((C(1,2)\^\,2-C(1,1)*C(2,2)))/p(3);

e(5)=sqrt((C(3,2)\^\,2-C(3,3)*C(2,2)))/p(1);

e(6)=sqrt((C(1,3)\^\,2-C(1,1)*C(3,3)))/p(2);

for j=1:6

    K=C-e(j)*P;

    r(j)=rank(K);

    l=mean(K(:,1));

    n=mean(K(:,2));

    q=mean(K(:,3));

    [x,y]=solve('H(1,1)*x\^\,2+H(2,2)*y\^\,2+2*H(1,2)*x*y+...\\    
    2*H(1,3)*x+2*H(2,3)*y=0','l*x+n*y+q=0');

    m1=eval(sqrt(x(1)\^\,2+y(1)\^\,2));

    m2=eval(sqrt(x(2)\^\,2+y(2)\^\,2));

    w(j)=max(m1,m2);

end

for j=1:6

    if isreal(w(j))\&\& w(j)$\le$ nor\&\&(r(j)==1)

        numr=w(j);

    end

end

%

\bibliography{mybib}

\end{document}